\documentclass[12pt]{amsart}
\usepackage{amsmath,amsthm,amssymb,color}

\def\pmatrix{\left(\begin{matrix}}
\def\endpmatrix{\end{matrix}\right)}

\def\Z{{\mathbb Z}}

\def\C{{\mathbb C}}

\def\Q{{\mathbb Q}}

\def\de{\delta}

\def\t{\theta}

\def\e{\varepsilon}

\def\e{\epsilon}

\def\A{{\mathcal A}}
\def\M{{\mathcal M}}
\def\J{{\mathcal J}}
\def\H{{\mathcal H}}

\def\tch#1#2{{\left[\begin{matrix}{#1}\\ {#2}\end{matrix}\right]}}
\def\stch#1#2{{\left[\begin{smallmatrix}{#1}\\ {#2}\end{smallmatrix}\right]}}

\def\Sp{\operatorname{Sp}(g,\Z)}

\def\dim{\operatorname{dim}}

\theoremstyle{plain}
\newtheorem{thm}{Theorem}
\newtheorem{lm}[thm]{Lemma}

\newtheorem{conj}{Conjecture}
\newtheorem{cor}[thm]{Corollary}

\theoremstyle{definition}
\newtheorem{df}{Definition}

\newtheorem{rem}[thm]{Remark}

\title {Affine stratification of $\A_4$.}
\author{Anant Atyam}

\address{Mathematics Department, Stony Brook University, Stony Brook, NY 11794-3651, USA.}
\date{\today} \email{anant@math.sunysb.edu}
\begin{document}
\begin{abstract}
We construct an affine stratification for $\A_4(\C)$ of length 6, in the sense of Roth and Vakil \cite{rova}, which gives us an upper bound of 6 for the cohomological dimension of $\A_4(\C)$. We conjecture that in general for arbitrary $g$ the cohomological dimension of $\A_g$ is equal to $g(g-1)/2$.
\end{abstract}

\maketitle

\section{Introduction}
It is well known that the coarse moduli space of principally polarized abelian varieties (ppav) $\A_g$ is not a projective variety over $\C$. Natural compactifications of $\A_g$ such as the Satake-Baily-Borel compactification $\A_g^{\rm Sat}$ \cite{satake} and various toroidal compactifications have been constructed, see for example \cite{namikawabook},\cite{fachbook}. The quasi projectivity of $\A_g$ raises  interesting questions, regarding its cohomological dimension and related invariants such as the affine covering number and affine stratification number. We address some questions regarding these invariants.

We now introduce some definitions \cite{rova} to make the comments above precise.
\begin{df}
The {\em affine covering number} (acn) of a scheme $X$ is one less than the smallest number of open affine sets required to cover $X$.
\end{df}

\begin{df}
The {\em coherent cohomological dimension} $\operatorname{cd}(X)$ of a variety $X$ is the smallest $i$ such that $H^j(X, \mathcal F) = 0$, for all $j>i$ and for any quasicoherent sheaf $\mathcal F$ on $X$.
\end{df}

Similarly

\begin{df}
The {\em constructible cohomological dimension} $\operatorname{ccd}(X)$ of a variety $X$ is the smallest $i$ such that $H^j(X, \mathcal F) = 0$ for all $j>i$ and for any constructible sheaf $\mathcal F$ on $X$.
\end{df}

\begin{df}
An {\em affine stratification} of a scheme $X$ is a finite decomposition $X = \sqcup_{k \in \mathbb Z \geq 0, i}Y_{k,i}$ into locally closed affine subschemes $Y_{k,i}$, where for each $Y_{k,i}$,
$$
  \overline Y_{k,i} \setminus Y_{k,i} \subseteq \cup_{k' > k, j}Y_{k',j}.
$$
The {\em length} of an affine stratification is the largest $k$ such that $\cup_j Y_{k,j}$ is non-empty.

The {\em affine stratification number} (asn) of a scheme $X$ is the minimum of the lengths of all possible affine stratifications of $X$.
\end{df}

The relationship between affine covering number, affine stratification number and cohomological dimension is given by the following theorem.

\begin{thm} \label{asnqcoh} \cite{rova}
The coherent cohomological dimension of a scheme $X$ is bounded above by the affine stratification number of the scheme $X$ which in turn is bounded above by the affine covering number of $X$, i.e., $\operatorname{cd}(X) \leq \operatorname{asn}(X) \leq \operatorname{acn}(X)$. Similarly the constructible cohomological dimension is bounded above by $\dim X + \operatorname{asn} X$, i.e. $\operatorname{ccd} X \leq \dim X + \operatorname{asn} X$.
\end{thm}

Our main result is the following statement, over $\C$.

\begin{thm}\label{thm:main}
There exists an affine stratification of $\A_4(\mathbb C)$ of length 6.
\end{thm}
Using theorem \ref{asnqcoh}, we get the following corollary:
\begin{cor}
The cohomological dimension $\operatorname{cd}(\A_4(\mathbb C)) \leq \operatorname{asn}(\A_4(\mathbb C)) \leq 6$, while the constructible cohomological dimension $\operatorname{ccd}(\A_4(\C)) \leq 16$.
\end{cor}

As an application of the study of the tautological ring of $\A_g$, we also prove the following proposition.

\begin{thm}\label{affine covering number}
The affine covering number of $\A_g$ for all $g$ is bounded below by $g(g-1)/2$, i.e., $\operatorname{acn}(\A_g) \geq g(g-1)/2$.
\end{thm}

\begin{rem}\label{remark}
Let $\A_g(\mathbb F_p)$ denote the moduli space of ppav over a finite field of characteristic $p$. Oort \cite{oortconj} showed that $\A_g(\mathbb F_p)$ has a complete subvariety $Y_p$ of codimension $g$ --- the locus of ppav the scheme of whose $p$-torsion points is supported at the origin (completeness follows from the fact that $\mathbb G_m$ always has non-trivial $p$-torsion, so that no semiabelic variety can lie in the closure). Since $Y_p$ is complete, we have $\operatorname{cd}(Y_p) = g(g-1)/2$, and thus since $Y_p$ is a closed subvariety of $\A_g(\mathbb F_p)$, it follows that $\operatorname{cd}(\A_g(\mathbb F_p)) \geq \operatorname{cd}(Y_p) = g(g-1)/2$. In particular we have $\operatorname{cd}(\A_4(\mathbb F_p))\ge 6$.
\end{rem}

\begin{rem}
Oort \cite{oortconj} had conjectured that $\A_g(\mathbb C)$ does not have a complete subvariety of codimension $g$, this conjecture was proven by Keel and Sadun in \cite{kesa}.
\end{rem}

Based on Remark \ref{remark} and  Theorem \ref{affine covering number} we make the following conjecture.

\begin{conj}\label{cohomological dim}
The cohomological dimension of $\A_g$ in all characteristics is equal to $g(g-1)/2$, i.e. $\operatorname{cd}(\A_g) = g(g-1)/2$.
\end{conj}
\begin{rem}
Fontanari and Pascolutti \cite{Fontanari:2011aa} constructed an affine open {\em cover} of the moduli space $\M_4$ of curves of genus 4. They did this by carefully constructing three modular forms whose zero loci have an empty intersection in $\M_4$. Two of the modular forms are combinatorially constructed so as to ensure the vanishing of two theta constants and a third modular form is constructed to avoid the hyperelliptic locus, which is defined in $\M_4$ by the condition that two theta constants vanish. We will use a similar approach to construct an affine stratification of $\A_4$, but after some point, our strategy for the construction will diverge from Fontanari and Pascolutti's approach. Unfortunately the strategy used by Fontanari and Pascolutti will not work for our purposes because of the inherent combinatorial difficulties in understanding vanishing loci of theta constants with characteristics.
\end{rem}

\subsection*{Acknowledgements}
I would like to thank my advisor Samuel Grushevsky for suggesting this problem to me. I am grateful to Klaus Hulek and Riccardo Salvati Manni for going through the pre-print meticulously and giving valuable comments and suggestions.

\section{Notation and Preliminaries}

We will exclusively work over $\C$ and we denote by $\H_g$ the Siegel upper half space of $g \times g$ complex symmetric matrices with positive definite imaginary part. There is an action of $\Sp$ on $\H_g$ given by
\begin{equation}\label{ABCD}
  \left(\begin{smallmatrix} A &B \\  C & D \end{smallmatrix}\right) \circ \tau = (A\tau + B)(C\tau + D)^{-1}.
\end{equation}
Here elements of $\Sp$ are written in block form and they consist of the integer matrices that preserve the skew symmetric non-degenerate bilinear form given by  $\left(\begin{smallmatrix}0& 1_g\\ -1_g&0
\end{smallmatrix}\right)$.

We recall that $\A_g$ has the structure of a quasi-projective variety (see \cite{igusabook}, \cite{grAgsurvey}) and can also be described as a group quotient $\A_g = \Sp\backslash\H_g$. A scalar {\em modular form} of weight $k$ with respect to a subgroup $\Gamma \subset \Sp$ is a holomorphic function $F: \H_g \to \C$ such that
\begin{equation}\label{modular}
  F((A\tau +B)(C\tau + D)^{-1}) = \det(C\tau + D)^kF(\tau).
\end{equation}
for any $\left(\begin{smallmatrix}A&B\\ C&D \end{smallmatrix}\right) \in \Gamma$.

We denote by $\t_m(z, \tau)$ the {\em theta function with characteristic} $m= \stch\e\de$ given by
\begin{equation}
  \t_m(z, \tau) := \sum_{n \in \Z^g}e^{\pi i((n+\e/2)^t\tau(n+\e/2) + 2(n+\e/2)^t(z+\de/2))}.
\end{equation}
Here $\stch\e\de \in (\Z\backslash2\Z)^{2g}$ and $z \in \C^g$. We shall say $m$ is an {\em even} characteristic if $  e(m) := \e^t\de = 0 \in\Z_2$.\label{parity}. We will denote $\mathcal E$ to be the set of all even characteristics.

It can easily be checked that $\t_m(z, \tau)$ is an even function of $z$ if and only if $m$ is an even characteristic. We shall call $\t_m(0, \tau)$ the {\em theta constant} with characteristic $m$. It follows from the above discussion that $\t_m(0, \tau)$ vanishes identically if $m$ is not an even characteristic. We shall denote $\mathcal E$ to be the set of all even theta characteristics in $(\Z\backslash2\Z)^8$. From now on, whenever we refer to theta constants, it is implicitly understood that we refer to theta constants with even characteristics.

It is known that theta constants are modular forms of weight $1/2$ with respect to a certain subgroup $\Gamma(4, 8) \subset \Sp$\cite{igusabook}. For our purposes we will also need the fact that $\t_m(0, \tau)^8$ is a modular form of weight 4 with respect to $\Gamma(2) \subset \Sp$. Here $\Gamma(2)$ is the level 2 congruence subgroup of $\Sp$ given by
\begin{equation}
\Gamma(2) := \left \{M \in \Sp \mid M \equiv 1_{2g} \mod 2  \right \}.
\end{equation}
The {\em level two cover} of $\A_g$ is defined as
$$
 \A_g(2) := \Gamma(2)\backslash\H_g.
$$
The Galois group of the cover $p:\A_g(2) \to \A_g$ is equal to $\operatorname{Sp}(g, \Z\backslash2\Z)$, and the above discussion means that eighth powers of theta constants are sections of an appropriate line bundle on $\A_g(2)$.

Under the action of the entire modular group $\Sp(g,\Z)$ theta constants transform as follows. We define the affine action of $\operatorname{Sp}(g, \Z\backslash2\Z)$ on theta characteristics by the formula \cite{igusabook}
\begin{equation}
\gamma *\tch\e\de = \begin{pmatrix}\label{trans}
D&C\\
B&A
\end{pmatrix}
\tch\e\de +\tch{\operatorname{diag}AB^t}{\operatorname{diag}CD^t}.
\end{equation}
for $\gamma = \left(\begin{smallmatrix}
A&B\\
C&D
\end{smallmatrix}\right) \in \Sp$, then in this notation, we have
\begin{equation}\label{automorphy}
\t_m(0, \gamma\circ\tau)^8 = \det(C\tau+D)^4\t_{\gamma*m}(0, \tau)^8.
\end{equation}
Thus a polynomial in theta constants defines a modular form for the entire $\Sp(g,\Z)$ if it is invariant under the affine action on characteristics

We shall denote the vector space of scalar modular forms of weight $k$ with respect to $\Sp$ by $M_{k}$. There is a natural compactification of $\A_g$ called the {\em Satake} compactification defined by
$$
  \A_g^{\rm Sat} := \operatorname{Proj}\bigoplus_{k \in \mathbb Z_{ \geq 0}}M_{k}.
$$
Thus using the definition, one sees that zero loci of scalar modular forms with respect to $\Sp$ are ample divisors in $\A_g^{\rm Sat}$. The modular forms that we will consider are all cusp forms i.e.~ they all vanish identically on the boundary $\A_g^{\rm Sat} \setminus \A_g$. Thus in particular if $D \subset \A_g^{\rm Sat}$ is a closed subscheme and if $F$ is a cusp form such that the zero locus $Z(F)$ does not contain any irreducible components of $D$ then $Z(F) \cap D$ is an ample divisor in $D$ and $D \setminus Z(F) \subset \A_g$ is affine.

We will denote $\A_g^{dec}$ to be the Zariski closed subset of $\A_g$ consisting of decomposable ppav, i.e.~ ppav that are isomorphic to a product of lower-dimensional ppav. We will denote $\A_g^{ind}:=\A_g\setminus\A_g^{dec}$ the Zariski open subset of $\A_g$ consisting of ppav that are not isomorphic to a product of lower dimensional ppav. We will denote $Hyp_g$ the locally closed subvariety of $\A_g$ consisting of the hyperelliptic jacobians.

\begin{rem}
It is worthwhile to recall that $Hyp_g$ is affine. This is seen to be the case as hyperelliptic curves of genus $g$ are double covers of $\mathbb P^1$ branched along $2g+2$ distinct points on $\mathbb P^1$. Thus a hyperelliptic curve is determined by a choice of $2g+2$ points on $\mathbb P^1$ up to an action of $\operatorname{Aut}(\mathbb P^1)$. Thus $Hyp_g$ can be naturally identified with a finite group quotient of $\M_{0, 2g+2}$, here $\M_{0,2g+2}$ denotes the moduli space of ordered $2g+2$ points on $\mathbb P^1$. The space $\M_{0,2g+2}$ is known to be affine, see for example \cite{hamobook}.
\end{rem}

From now on we will abuse notation and denote $\A_k \times \A_{g-k}$ the closed subset of $\A_g$ consisting of ppav that are isomorphic to a product of a $k$-dimensional ppav and a  ($g-k$)-dimensional ppav. Strictly speaking this is an abuse of notation as we could have used $\A_{g-k} \times \A_k$ to denote this closed subset as well. Thus for example when we write $\A_2 \times \A_2 \subset \A_4$ we really mean $\operatorname{Sym}^2(\A_2) \subset \A_4$.

We will now state a result of Igusa and Salvati Manni on the orbits of tuples of theta characteristics under the action of $\operatorname{Sp}(g, \Z\backslash2\Z)$ given by \eqref{trans}, which we will use in our construction of modular forms.

\begin{thm}[\cite{igusabook}, see also \cite{smlevel2}]\label{tuplesorbits}Given two p-tuples of {\em even} theta characteristics $I = (m_1, m_2, ... m_p)$ and $J = (n_1, n_2.... n_p)$, they are in the same orbit of the action of the deck group $\operatorname{Sp}(2g,\Z\backslash2\Z)$ on the set of p-tuples of theta characteristics, if and only if the following conditions hold:

1. Linear relations among even number of terms in $I$ translate to corresponding linear relations in $J$, i.e., $m_{i_1} + m_{i_2} .... + m_{i_{2l}} = 0$  if and only if $n_{i_1} + n_{i_2} + n_{i_3} .... n_{i_{2l}} = 0$.

2. $ e(m_i+m_j+m_k) = e(n_i+n_j+n_k)$  $\forall i, j, k$.
\end{thm}

For a given $k$, let $n_k$ be the number of theta characteristics in $\mathcal E$ of the form $\tch\e\de = \begin{bmatrix}\e_{1} & \e_{2}\\ \de_{1} & \de_2 \end{bmatrix}$\label{oddk} with $\tch {\e_{1}} {\de_{1}} \in (\Z\backslash2\Z)^{2k}_{odd}$ and $\tch {\e_2} {\de_2} \in (\Z\backslash2\Z)^{2(g-k)}_{odd}$, i.e. $\e_1^{t}\de_1 = \e_2^{t}\de_2 = 1 \in \Z\backslash2\Z$.

Let $I \in \mathcal E^{n_k}$ be a $n_k$ tuple all of whose characteristics are of the form described above. Then we have the following theorem.

\begin{thm}[\cite{smlevel2}]
Given $\tau \in \H_g$ then $[\tau]$ corresponds to a ppav that is a product of a ppav of dimension $k$ and $g-k$ if and only iff $\exists J = (m_1, m_2 \ldots m_{n_k}) \in \mathcal E^{n_k}$ in the $\operatorname{Sp}(g, \Z\backslash2\Z)$ orbit of $I$ with $\t_{m_1}(0, \tau) = \t_{m_2}(0, \tau) ..... \t_{m_{n_k}}(0, \tau) = 0$.
\end{thm}

Consider the universal family $\pi: \mathcal X_g \to \A_g$. The Hodge bundle $\mathbb E$ is the locally free sheaf on $\A_g$ defined by the sheaf of relative differentials, $E:= \pi_*(\Omega^1_{\mathcal X_g /\A_g})$ and let $\lambda_i$ be the $i$'th Chern class of this sheaf. It has been shown by van der Geer \cite{vdgeercycles} and Esnault and Viehweg \cite{esvi} that $\lambda_1^{g(g-1)/2} \neq 0 \in CH^*_{\Q}(\A_g)$. It is also known that $\operatorname{Pic}_\Q(\A_g)$ is one dimensional for $g \geq 2$, \cite{borel}. This implies that there does not exist a collection of $g(g-1)/2$ divisors in $\A_g$ with an empty intersection.
We claim that this implies that $\operatorname{acn}(\A_g) \geq g(g-1)/2$, because of the following lemma.

For $g \geq 2$,

\begin{lm}If $U$ is an open affine subset of $\A_g$ then $U \subset \A_g \setminus D$, where $D$ is a divisor.
\begin{proof}Let $i: U \to \A_g$ be the inclusion map. We know that $\operatorname{Pic}_{\Q}(\A_g) = \Q\lambda_1$\cite{borel}. We also know that since $\A_g$ is normal, $\operatorname{Pic}_{\Q}(U) = 0$. The only way for $i^*:Pic_{\Q}(\A_g) \to Pic_{\Q}(U)$ to be the zero map is if the $\lambda_1$ class goes to zero and that can happen only if $U$ is contained in the complement of the divisor, as removing a higher codimension locus of $\A_g$ does not affect the image of $\operatorname{Pic}_{\Q}(\A_g)$.
\end{proof}
\end{lm}
We have thus proven
\begin{thm}
The affine covering number $\operatorname{acn}(\A_g) \geq g(g-1)/2$.
\end{thm}

\section{Proof of main theorem}
The proof of Theorem 2 is divided into two steps. The first step consists of writing down three modular forms whose zero loci in $\A_4$ intersect in a Zariski closed subset of $\A_4$ of codimension 3, thus giving us three affine strata. This construction is similar to that of Fontanari and Pascolutti \cite{Fontanari:2011aa}. The second step consists of constructing the remaining 4 strata by identifying natural codimension one subvarieties of closures of successive strata.

\begin{proof}
\textbf {Step 1(Strata induced by zero loci of modular forms)}:

\textit{1st Stratum}: Consider the Schottky form defined by the following formula, using the notation in \cite{Fontanari:2011aa} \begin{equation}\label{Schottky}
 F_T := 16\sum_{m \in \mathcal E}\theta_m^{16} - \sum_{m \in \mathcal E}(\theta_m^8)^2.
\end{equation}
Igusa \cite{igusachristoffel} showed that the zero locus $Z(F_T)$ of $F_T$ is equal to $\J_4 \subset \A_4$. We define our first stratum to be
$$
  X_0 := \A_4 \setminus Z(F_T).
$$
\textit{2nd Stratum}: Consider the Theta null form on $\A_4$
\begin{equation}
 \Theta_{null} = \prod_{m \in \mathcal E}\theta_m.
\end{equation}
From Riemann's theta singularity theorem it easily follows that the modular form $\Theta_{null}$ does not vanish identically on $\M_g$ in any genus. We thus define our second stratum to be
$$
  X_1 := Z(F_T) \setminus Z(\Theta_{null}).
$$
This stratum corresponds to Jacobians of canonical curves of genus four that are given by a complete intersection of a smooth cubic and a smooth quadric. \cite{dolgachevbook}

\begin{rem}

The group action of $\operatorname{Sp}(g, \Z\backslash2\Z)$ on theta characteristics becomes important from the third stratum onwards. The reason for its importance is the following:
if $V_1 = Z(\theta_{m_1}^8, \theta_{m_2}^8) \subset \A_g(2)$ and if $V_2 = Z(\theta_{n_1}^8, \theta_{n_2}^8)$, then (\ref{automorphy}) and Theorem \ref{tuplesorbits} tell us that $V_1$ and $V_2$ are equivalent to each other under the action of the Galois group of the cover $\A_g(2) \to \A_g$.

In general, from the transformation formula for theta constant \eqref{automorphy} of theta constants, it follows that if $m_i$'s and $n_i$'s satisfy the conditions of Theorem \ref{tuplesorbits}, then the Galois group of this level 2 cover conjugates $Z(\theta_{m_1}^8, \theta_{m_2}^8, .... \theta_{m_k}^8)\subset\A_g(2)$ to $Z(\theta_{n_1}^8, \theta_{n_2}^8 .... \theta_{n_k}^8)\subset\A_g(2)$ .  We should mention that R. Varley constructed an example of a set of 10 theta constants vanishing on a level cover of $\A_4$ that describes a single ppav, which corresponds to the Segre cubic threefold with an even point of order two; i.e., a cubic threefold with 10 nodes \cite{varley}. The correspondence between cubic threefolds and $\A_4$ arises due to the work of Donagi and Smith \cite{dosm}, further explored by Izadi \cite{izadiA4}, which establishes a birational map between $\A_4$ and a level cover of the moduli of cubic three folds $\mathcal C$. This is in contrast to the situation where a different set of 10 theta constants i.e. one lying in a different orbit of the $\Sp$ action in (\ref{trans}) describes the hyperlliptic locus.\end{rem}

\textit{3rd Stratum}: Consider the modular form $F_1$ given by
\begin{equation}
F_1 = \sum_{m \in \mathcal E}\frac{\Theta_{null}^8}{\theta_m^8}.
\end{equation}

The Zariski closed subset given by $Z(\Theta_{null}) \cap Z(F_1) \cap Z(F_T)$ can be interpreted as the closed subset of $\J_4$ where some two theta constants vanish. Igusa \cite{igusagen4} proved that a curve $[C] \in \M_4 \subset \A_4$ is hyperelliptic if and only if two theta constants vanish on the Jacobian of $C$. On the other hand Theorem 7 and Remark 8 tell us that $\A_4^{dec} \subset Z(F_T) \cap Z(\Theta_{null}) \cap Z(F_1)$. Hence we see $Z(F_T) \cap Z(\Theta_{null}) \cap Z(F_1) = Hyp_4 \cup \A_4^{dec}$. We thus define our third stratum to be
$$
  X_2 := (Z(F_T) \cap Z(\Theta_{null})) \setminus Z(F_1).
$$

This stratum corresponds to Jacobians of canonical curves of genus four that are a complete intersection of a smooth cubic and a singular quadric of rank three. \cite{dolgachevbook}

\begin{rem}
The two modular forms given by $F_1$ and $\Theta_{null}$ were used by Fontanari and Pascolutti in proving Looijenga's conjecture for $\M_4$ \cite{Fontanari:2011aa}. They then use another modular form $F_H$ again constructed as a sum of products of theta constants. They then go on to show that this modular form does not intersect $Hyp_4 \subset \M_4$. { \em We will not be able to use $F_H$ to construct our 4th stratum as $Z(F_T, F_1, \Theta_{null}, F_H) = \A_4^{dec}$}, since the codimension of $\A_4^{dec} \subset \A_4$ equals to 3, we see that the stratification using $F_H$ is not optimal as we will be using four forms to describe a geometric locus of codimension 3.
\end{rem}

\textbf{Step 2 (Affine stratification constructed by identifying appropriate closed subvarieties)}

\textit{4th Stratum}: Consider the boundary of the previous stratum, which is $\overline{X_2} \setminus X_2
\ =\ Hyp_4 \cup \A_4^{dec}$. We know that $Hyp_g$ is affine for all $g$. The boundary of $\overline{Hyp_4}$  is equal to $(\A_1 \times \overline{Hyp_3})\ \cup\ (\A_2 \times \A_2)$. It is known classically that in $\A_3$, jacobians of hyperelliptic curves are characterized by the vanishing of a single theta constant, that is to say $\overline{Hyp_3}=Z(\Theta_{null})$. Since moreover $\A_1$ is affine, the product $\A_1 \times ({\A_3 \setminus \overline{Hyp_3}})$ is affine, as a product of affine varieties. We thus define our fourth stratum to be
$$
  X_3 := Hyp_4\ \sqcup\ \left(\A_1 \times {\A_3 \setminus \overline{Hyp_3}})\right.
$$
Since $X_3$ is a disjoint union of two affine varieties, it is itself affine.

\textit{5th Stratum}: Consider the boundary of the previous stratum
$$
  \overline{X_3} \setminus X_3\ =\ (\A_2 \times \A_2) \cup (\A_1 \times \overline{Hyp_3}).
$$
Now recall that that any indecomposable abelian threefold is a Jacobian of a smooth curve; it thus follows that $\overline{Hyp_3} \setminus Hyp_3 = \A_1 \times  \overline{Hyp_2}$, while $\A_2 \times \A_2 = \M_2^{ct} \times \M_2^{ct}$. Thus by remark 5 we see that
$$
 \left(\A_2 \times \A_2\right) \setminus \left(\A_1 \times \A_1 \times \A_2\right) = \M_2 \times \M_2.
$$
We know $\M_2$ to be affine as $\M_2 = Hyp_2$. On the other hand, we have
$$
 \left(\A_2 \times \A_2\right) \cap \left(A_1 \times \overline{Hyp_3}\right) = \A_1 \times \A_1 \times \A_2.
$$
We thus define our fifth stratum to be
$$
  X_4 := \left(\A_1 \times Hyp_3\right) \sqcup \left(\M_2 \times \M_2\right).
$$
This stratum is seen to be affine as $\A_1 = Hyp_1$ is affine, as are $Hyp_3$ and $\M_2$.

\textit{6th Stratum}: Consider the boundary of the previous stratum $\overline {X_4} \setminus X_4 = \A_1 \times \A_1 \times \A_2$. Since $\A_1$ is known to be affine we just need to stratify $\A_2$. But $\A_2 = \M_2^{ct}$, and $\M_2^{ct} = \M_2 \sqcup \A_1 \times \A_1$. Thus we define our sixth stratum to be
$$
  X_5 := \left(\A_1 \times \A_1 \times \A_2\right) \setminus \left(\A_1 \times \A_1 \times \A_1 \times \A_1\right).
$$
This stratum is affine as it simply equal to $\A_1 \times \A_1 \times \M_2$.

\textit{7th Stratum}: We finally define our seventh stratum to be
$$
  X_6:= \A_1 \times \A_1 \times \A_1 \times \A_1.
$$
We then see that $\A_1$ is affine as $\A_1 = Hyp_1$, and the theorem is thus proven.
\end{proof}

\begin{rem}
The next question that is natural to consider is if there exists an affine stratification of $\A_5$. There are a few difficulties in approaching this problem. For one, the geometry of $\M_5$ in $\A_5$ is not easily accessible through the machinery of modular forms as $\operatorname{codim}(\M_5 \subset \A_5) = 3$ and it is not known if the Schottky locus is a set theoretic complete intersection in $\A_5$.
\end{rem}

\bibliographystyle{alpha}
\bibliography{sam_biblio}

\end{document}